\documentclass{amsart} 
 \newtheorem{thm}{Theorem}[section]
 \newtheorem{cor}[thm]{Corollary}
 \newtheorem{lem}[thm]{Lemma}
 
 \theoremstyle{definition}
 
 \theoremstyle{remark}

 \numberwithin{equation}{section}

\newcommand{\scal}[1]{\langle#1\rangle}
\newcommand{\norm}[1]{\left\Vert#1\right\Vert}
\newcommand{\abs}[1]{\left\vert#1\right\vert}
\newcommand{\trace}{\operatorname{tr}}
\newcommand{\C}{\mathbb C}
\newcommand{\R}{\mathbb R}
\newcommand{\dist}{\operatorname{dist}}
\renewcommand{\Re}{\operatorname{Re}}

\begin{document}
%
%
%
%
%
%
%
%
%
\title[Numerical Ranges of Quadratic Operators]
 {On Generalized Numerical Ranges\\ of Quadratic Operators}
\author[Rodman]{Leiba Rodman}

\address{%
Department of Mathematics \\
College of William and Mary \\
Williamsburg, VA 23185\\
USA} \email{lxrodm@math.wm.edu}

\thanks{The research of both authors was partially supported by NSF grant DMS-0456625.}
\author[Spitkovsky]{Ilya M. Spitkovsky}
\address{%
Department of Mathematics \\
College of William and Mary \\
Williamsburg, VA 23185\\ USA}\email{ilya@math.wm.edu}


\date{}

\begin{abstract}
It is shown that the result of Tso-Wu on the elliptical shape of
the numerical range of quadratic operators holds also for the
essential numerical range. The latter is described quantitatively,
and based on that sufficient conditions are established under
which the $c$-numerical range also is an ellipse. Several examples
are considered, including singular integral operators with the
Cauchy kernel and composition operators.
\end{abstract}

\maketitle
\section{Introduction}
Let $A$ be a bounded linear operator acting on a complex Hilbert
space $\mathcal H$. Recall that the {\em numerical range} $W(A)$
of $A$ is defined as \[ W(A) =\{ \scal{Ax,x}\colon x\in{\mathcal
H}, \norm{x}=1 \} .\] If $c$ is a $k$-tuple of non-zero (in
general, complex) numbers $c_1,\ldots,c_k$, then the $c$-{\em
numerical range} of $A$ is \[ W_c(A)=\left\{ \sum_{j=1}^k c_j
\scal{Ax_j,x_j}\colon \{ x_j\}_{j=1}^{k}\text{ is an orthonormal
subset of } {\mathcal H}\right\}.\] Of course, if $c$ consists of
just one number $c_1=1$, $W_c(A)$ is nothing but the regular
numerical range of $A$. Also, for $c_1=\ldots=c_k=1$, the
$c$-numerical range $W_c(A)$ turns into $W_k(A)$ -- the so called
$k$-{\em numerical range} \footnote{We realize that there is a
slight abuse of notation here, but both $W_c(A)$ and $W_k(A)$ are
rather standard, and the meaning is usually clear from the
content.}  introduced by Halmos, see \cite{Ha}. Finally, the {\em
essential numerical range} introduced in \cite{StWil68} can be
defined \cite{FilStWil72} as
\begin{equation}\label{wess} W_{\rm ess}(A)=\bigcap {\rm cl}\,
W(A+K),\end{equation} where the intersection is taken over all
compact on $\mathcal H$ operators $K$, and the symbol cl denotes
the topological closure. Considering $W_c(A)$ or $W_{\rm ess}(A)$,
we will implicitly suppose that $\dim {\mathcal H}\geq k$ or that
$\mathcal H$ is infinite dimensional, respectively.

There are several monographs devoted to the numerical range and
its various generalizations (including those mentioned above), see
for example \cite{BoDu73,GusRa}. We mention here only the results
which are of direct relevance to the subject of this paper.

From the definitions it is clear that all three sets are unitarily
invariant:

\begin{equation}\label{uninv} W(U^*AU)= W(A),\ W_c(U^*AU)= W_c(A),
\ W_{\rm ess}(U^*AU)= W_{\rm ess}(A)\end{equation} for any unitary
operator $U$ on $\mathcal H$. Also, they behave in a nice and
predictable way under affine transformations of $A$:
\begin{equation}\label{aff} W(\alpha A+\beta I) =\alpha W(A)+\beta, \ W_{\rm
ess}(\alpha A+\beta I) =\alpha W_{\rm ess}(A)+\beta,
\end{equation} and
\begin{equation}\label{affc}
W_c(\alpha A+\beta I) =\alpha W_c(A)+\beta\sum_{j=1}^k c_j
\end{equation} for any $\alpha,\beta\in\C$.

It is a classical result (known as the Hausdorff-Toeplitz theorem)
that the set $W(A)$ is convex. Clearly, $W_{\rm ess}(A)$ is
therefore convex as well. The $c$-numerical range is convex if all
$c_j$ lie on the same line passing through the origin but not in
general \cite{West75}. In what follows, we suppose that $c_j$
satisfy the above mentioned condition. Moreover, since \[
W_c(\alpha A)= W_{\alpha c}(A), \quad \alpha\in\C , \] we then may
(and will) without loss of generality suppose that all $c_j$ are
real. We will also arrange them in the non-increasing order: \[
c_1\geq c_2\ldots\geq c_k,\] since permutations of $c_j$ leave
$W_c(A)$ invariant.

When $\dim{\mathcal H}=2$, the numerical range of $A$ is the
closed (as is always the case in finite dimensional setting)
elliptical disc with the foci at the eigenvalues $\lambda_1,
\lambda_2$ of $A$ and the minor axis $\sqrt{\trace
(A^*A)-\abs{\lambda_1}^2-\abs{\lambda_2}^2}$ (the elliptic range
theorem, see, e.g., \cite[Section 1.1]{GusRa}). According to the
Cayley - Hamilton theorem, $A$ in this setting satisfies the
equation
\begin{equation}\label{quad} A^2-2\mu A-\nu I=0 \end{equation}
with \[ \mu=(\lambda_1+\lambda_2)/2,\quad
\nu=-\lambda_1\lambda_2.\] For arbitrary $\mathcal H$, operators
$A$ satisfying (\ref{quad}) with some $\mu,\nu\in\C$ are called
{\em quadratic operators}.

Rather recently, Tso and Wu showed that $W(A)$ is an elliptical
disc (open or closed) for any quadratic operator $A$, independent
of the dimension of $\mathcal H$ \cite{TsoWu}.

In this paper, we continue considering the (generalized) numerical
ranges of quadratic operators. We start by stating Tso-Wu's result
and outlining its proof (different from one presented in
\cite{TsoWu}), in order to show how it can be modified to prove
ellipticity of the {\sl essential} numerical ranges of quadratic
operators. We then use the combination of the two statements to
derive some sufficient conditions for the $c$-numerical range to
also have an elliptical shape. This is all done in Section~1.
Section~2 is devoted to concrete implementations of these results.

\section{Main results}

\subsection{Classical numerical range}
We begin with the Tso-Wu result.
\begin{thm}Let the operator $A$ satisfy equation {\em (\ref{quad})}.
Then $W(A)$ is the elliptical disc with the foci
$\lambda_{1,2}=\mu\pm\sqrt{\mu^2+\nu}$ and the major/minor axis of
the length \begin{equation}\label{axes}
s\pm\abs{\mu^2+\nu}s^{-1}.\end{equation} Here $s=\norm{A-\mu I}$,
and the set $W(A)$ is closed when the norm $\norm{A-\mu I}$ is
attained and open otherwise.\label{th:w}\end{thm} \begin{proof} As
in \cite[Theorem 1.1]{TsoWu}, observe first that (\ref{quad})
guarantees unitary similarity of $A$ to an operator of the form
\begin{equation}\label{canA} \lambda_1 I\oplus \lambda_2 I\oplus
\left[\begin{matrix} \lambda_1 I & 2X \\ 0 & \lambda_2 I
\end{matrix}\right] \end{equation} acting on
${\mathcal H}_1\oplus{\mathcal H}_2\oplus ({\mathcal
H}_3\oplus{\mathcal H}_3)$, where $\dim{\mathcal H}_j\, (\geq 0)$
are defined by $A$ uniquely, and $X$ is a positive definite
operator on ${\mathcal H}_3$. According to the first of properties
(\ref{uninv}), we may suppose that $A$ itself is of the form
(\ref{canA}).

Using the first of formulas (\ref{aff}) we may further suppose
that $\mu=0$ and $\nu\geq 0$; in other words, that in (\ref{canA})
\begin{equation}\label{la}\lambda_1=-\lambda_2 :=\lambda\geq
0,\quad \lambda^2=\nu.\end{equation}

The case ${\mathcal H}_3=\{0\}$ corresponds to the normal operator
$A$ when $W(A)$ is the closed line segment connecting $\lambda_1$
and $\lambda_2$. This is in agreement with formula (\ref{axes})
when $\nu\neq 0$, since in this case $s=\sqrt{\nu}$ is attained,
and $s-\sqrt{\nu}s^{-1}=0$. In the trivial case $s=0$ (when the
operator $A$ is scalar and $W(A)$ degenerates into a single point)
formula (\ref{axes}) formally speaking is not valid since $s^{-1}$
is not defined. However, the relation between $s$ and $\nu$
justifies the convention $\nu s^{-1}=0$ in this case.

In the non-trivial case $\dim{\mathcal H}_3>0$ our argument is
different from that in \cite{TsoWu}. Namely, we will make use of
the fact that the (directed) distance from the origin to the
support line ${\ell}_\theta$ with the slope $\theta$ of $W(A)$ is
the maximal point $\omega_\theta$ of the spectrum of $\Re
(ie^{-i\theta}A)$. Moreover, $\ell_\theta$ actually contains
points of $W(A)$ if and only if $\omega_\theta$ belongs to the
point spectrum of $\Re (ie^{-i\theta}A)$.

For $A$ of the form (\ref{canA}) with $\lambda_j$ as in
(\ref{la}), \[ \Re (ie^{-i\theta}A) = (\lambda\sin\theta) I\oplus
(-\lambda\sin\theta) I\oplus\left[\begin{matrix}
(\lambda\sin\theta) I & ie^{-i\theta}X
\\ -ie^{i\theta}X & (-\lambda\sin\theta) I
\end{matrix}\right]  .\] Thus,
\begin{multline}\label{re} \Re (ie^{-i\theta}A)-\omega I = \\ (\lambda\sin\theta-\omega)
I\oplus (-\lambda\sin\theta-\omega) I\oplus\left[\begin{matrix}
(\lambda\sin\theta-\omega) I & ie^{-i\theta}X
\\ -ie^{i\theta}X & -(\lambda\sin\theta+\omega) I
\end{matrix}\right]  .\end{multline} For any $\omega\neq
\lambda\sin\theta$, the last direct summand in (\ref{re}) can be
rewritten as \begin{equation}\label{id} 
\left[\begin{matrix} I & 0\\ 0 &
(\lambda\sin\theta-\omega)^{-1}I\end{matrix}\right]
\left[\begin{matrix} I & 0 \\ -ie^{i\theta}X &
I\end{matrix}\right] \left[\begin{matrix}
(\lambda\sin\theta-\omega) I & ie^{-i\theta}X \\ 0 &
(\omega^2-\lambda^2\sin^2\theta) I-X^2
\end{matrix}\right] .
\end{equation} Therefore, $\omega_\theta=\sqrt{\lambda^2\sin^2\theta+\norm{X}^2}$
is the rightmost point of the spectrum of $\Re (ie^{-i\theta}A)$.
In other words, the support lines of $W(A)$ are the same as of the
numerical range of the $2\times 2$ matrix
\[ \left[\begin{matrix}\lambda & 2\norm{X}\\ 0 &
-\lambda\end{matrix}\right] .\] The description of $W(A)$ as the
elliptical disc with the foci and axes as given in the statement
of the theorem follows from here and the elliptic range theorem.

Moreover, $\omega_\theta$ is an eigenvalue of $\Re
(ie^{-i\theta}A)$ if and only if the norm of $X$ (or equivalently,
of $A$ itself) is attained, so that this either happens for all
$\theta$ or for none of them. In the former case, every support
line of $W(A)$ must contain at least one of its points, and the
elliptical disc $W(A)$ is closed. In the latter case, the support
lines are disjoint with $W(A)$, so that it is open.
\end{proof}
{\bf Remark.} Formula (\ref{axes}) is formally different from the
result of \cite[Theorem 2.1]{TsoWu}, where the lengths of the axes
of $W(A)$ are given in terms of $\norm{A-\lambda_1 I}$, not
$\norm{A-\mu I}$. The two operators coincide when $\mu^2+\nu=0$.
If this is not the case, the relation between their norms follows
from the general property \[
\norm{P}=\frac{1}{2}(\norm{S}+\norm{S}^{-1}) \] of any projection
$P$ and associated with it involution $S=2P-I$ (see \cite{Spit80})
applied to $P=(A-\lambda_1 I)/(\lambda_2-\lambda_1)$ and $S=(A-\mu
I)/\sqrt{\mu^2+\nu}$.

As a matter of fact, the relation between $A$ and involution
operators shows that $A$ can be represented as a (rather simple)
function of two orthogonal projections. This observation allows to
describe the spectra and norms of all operators involved in the
proof of Theorem~\ref{th:w} straightforwardly, using the machinery
developed in \cite{Spitkov}. We chose an independent exposition,
in the interests of self containment.

\subsection{Essential numerical range}
If $A$ satisfies (\ref{quad}) and one of its
eigenvalues (say $\lambda_1$) has finite multiplicity, then in
representation (\ref{canA}) the spaces ${\mathcal H}_1$ and
${\mathcal H}_3$ are finite dimensional. Thus, $A$ differs from
$\lambda_2 I$ by a compact summand, and $W_{\rm ess}(A)$ is a
single point. Let us exclude this trivial situation, that is,
suppose that $\sigma_{\rm
ess}(A)=\sigma(A)=\{\lambda_1,\lambda_2\}$.

From (\ref{wess}) it is clear that the support lines
$\ell_\theta^{\rm ess}$ with the slope $\theta$ are at the
distance $\omega_\theta^{\rm ess}$ from the origin. Here
$\omega_\theta^{\rm ess}$ is the maximal point of the essential
spectrum of $\Re(ie^{-i\theta}A)$. This observation allows to
repeat the statement and the proof of Theorem~\ref{th:w} almost
literally, inserting the word ``essential'' where appropriate (of
course, the last paragraph of the proof becomes irrelevant since
the essential numerical range is always closed). We arrive at the
following statement.

\begin{thm}\label{th:wess}Let the operator $A$ satisfy equation {\em
(\ref{quad})}, with both eigenvalues
$\lambda_{1,2}=\mu\pm\sqrt{\mu^2+\nu}$ having infinite
multiplicity. Then $W_{\rm ess}(A)$ is the closed elliptical disc
with the foci $\lambda_{1,2}$ and the major/minor axis of the
length $s_0\pm\abs{\mu^2+\nu}s_0^{-1}$, where $s_0$ is the
essential norm of $A-\mu I$.\end{thm} In the trivial case $s_0=0$
(when $A$ differs from $\mu I$ by a compact summand, so that
necessarily $\mu^2+\nu=0$) we by convention set
$\abs{\mu^2+\nu}s_0^{-1}=0$. This agrees with the fact that
$W_{\rm ess}(A)$ then degenerates into a singleton $\mu$.

\begin{cor}\label{cor:clos}Let the operator $A$ satisfying
{\em (\ref{quad})} be such that \begin{equation}\label{ineqnor}
\norm{A-\mu I}>\norm{A-\mu I}_{\rm ess}.\end{equation} Then the
elliptical disc $W(A)$ is closed. \end{cor}
\begin{proof}
Indeed, (\ref{ineqnor}) holds if and only if $\norm{X}_{\rm
ess}<\norm{X}$ for $X$ from (\ref{canA}). Being positive definite,
the operator $X$ then has $\norm{X}$ as its eigenvalue. In other
words, the norm of $X$ (and therefore of $A-\mu I$) is attained.
It remains to invoke the last statement of
Theorem~\ref{th:w}.\end{proof}

\subsection{\boldmath{$c$}-numerical range}
The behavior of $W_c(A)$, even for quadratic operators, is more
complicated; see \cite{ChienTsoWu} for some observations on the
$k$-numerical range. With no additional assumptions on $A$, we
give only a rather weak estimate. In what follows, it is
convenient to use the notation $\norm{c}=\sum_{j=1}^k \abs{c_j}$.

\begin{lem}\label{l:wc}Let $A$ be as in Theorem~{\em
\ref{th:wess}}. Denote by $s$ and $s_0$ the norm and essential
norm of $A-\mu I$ respectively, and by $E$ and $E_0$ two
elliptical discs with the foci at $\mu\sum_{j=1}^k
c_j\pm\sqrt{\mu^2+\nu}\norm{c}$, the first -- closed, with the
axes $(s\pm\abs{\mu^2+\nu}s^{-1})\norm{c}$ and the second -- open,
with the axes $(s_0\pm\abs{\mu^2+\nu}s_0^{-1})\norm{c}$. Then
$W_c(A)$ contains $E_0$ and is contained in $E$.
\end{lem}
\begin{proof}Using (\ref{affc}) we may, as in the proof of
Theorem~\ref{th:w}, without loss of generality suppose that
$\mu=0$, $\nu\geq 0$. Since all the sets $E$, $E_0$ and $W_c(A)$
are convex, we need only to show that the support line to $W_c(A)$
in any direction lies between the respective support lines to
$E_0$ and $E$. In other words, the quantity
\begin{equation}\label{sup} \sup \left\{
\sum_{j=1}^kc_j\Re\scal{ie^{-i\theta}Ax_j,x_j}\colon \{
x_j\}_{j=1}^k\text{ is orthonormal}\right\}\end{equation} must lie
between $$\norm{c}\sqrt{\nu\sin^2\theta+\norm{X}_{\rm ess}^2}
\quad \mbox{and} \quad \norm{c}\sqrt{\nu\sin^2\theta+\norm{X}^2}$$
with $X$ given by (\ref{canA}). But this is indeed so, because
(\ref{id}) implies that the spectrum and the essential spectrum of
$\Re(ie^{-i\theta}A)$ have the endpoints
$\pm\sqrt{\nu\sin^2\theta+\norm{X}^2}$ and
$\pm\sqrt{\nu\sin^2\theta+\norm{X}_{\rm ess}^2}$, respectively.
\end{proof}

An interesting situation occurs when the norm of $A-\mu I$
coincides with its essential norm (equivalently,
$\norm{X}=\norm{X}_{\rm ess}$ for $X$ from (\ref{canA})), so that
$E$ is simply the closure of $E_0$. To state the explicit result,
denote by $m_\pm$ the number of positive/negative coefficients
$c_j$ and let $m=\max\{ m_+,m_-\}$.

\begin{thm}\label{th:wc}Let $A$ be as in Theorem~{\em
\ref{th:wess}}, and on top of that \begin{equation}\label{eqnor}
\norm{A-\mu I}=\norm{A-\mu I}_{\rm ess}. \end{equation} Define $E$
and $E_0$ as in Lemma $\ref{l:wc}$. Then $W_c(A)$ coincides with
$E$ if the norm of $A-\mu I$ is attained on the subspace of the
dimension at least $m$, and with $E_0$ otherwise.
\end{thm}
\begin{proof}Consider first a simpler case when in(\ref{canA})
$\dim{\mathcal H}_3<\infty$. Then due to (\ref{eqnor}) ${\mathcal
H}_3=\{0\}$, so that the operator $A$ is normal. The norm
$\abs{\mu^2+\nu}^{1/2}$ of $A-\mu I$ is attained on infinite
dimensional subspaces ${\mathcal H}_1$ and ${\mathcal H}_2$, and
$W_c(A)$ is the closed line segment connecting the points
$\mu\sum_{j=1}^k c_j+\sqrt{\mu^2+\nu}\norm{c}$ and
$\mu\sum_{j=1}^k c_j-\sqrt{\mu^2+\nu}\norm{c}$. This segment
apparently coincides with $E$.

Let now ${\mathcal H}_3$ be infinitely dimensional. From
Lemma~\ref{l:wc} it follows that $W_c(A)$ lies between $E$ and its
interior $E_0$, so that the only question is which points of the
boundary of $E$ belong to $W_c(A)$. It follows from (\ref{id})
that the minimal and maximal points of the spectrum of $\Re
(ie^{-i\theta}A)$ have the same multiplicity as its eigenvalues,
this multiplicity does not depend on $\theta$ and coincides in
fact with the dimension $d$ ($\geq 0$) of the subspace on which
the norm of $X$ is attained. From (\ref{canA}) under conditions
(\ref{la}) it follows that the norm of $A-\mu I$ is attained on a
$d$-dimensional subspace as well.

On the other hand, the supremum in (\ref{sup}) is attained if and
only if this multiplicity is at least $m$. Thus, the boundary of
$E$ belongs to $W_c(A)$ if $d\geq m$ and is disjoint with $W_c(A)$
otherwise.
\end{proof}

\section{Examples}
We consider here several concrete examples illustrating the above
stated abstract results. All the operators $A$ involved happen to
be involutions which corresponds to the choice $\mu=0$, $\nu=1$ in
(\ref{quad}). According to Theorems~\ref{th:w} and \ref{th:wess},
the major/minor axes of the elliptical discs $W(A)$ and $W_{\rm
ess}(A)$ then have the lengths \begin{equation}\label{invaxes}
\norm{A}\pm\norm{A}^{-1}\text{ and } \norm{A}_{\rm
ess}\pm\norm{A}_{\rm ess}^{-1}, \end{equation} respectively.

\subsection{Singular integral operators on closed curves}
Let $\Gamma$ be the union of finitely many simple Jordan
rectifiable curves. Suppose that the number of its points of
self-intersection is finite, and that $\Gamma$ partitions the
extended complex plane $\dot{\C}=\C\cup\{\infty\}$ into two open
disjoint (not necessarily connected) sets $D^+$ and $D^-$.
Moreover, we suppose that $\Gamma$ is the common boundary of $D^+$
and $D^-$, and that it is oriented in such a way that the points
of $D^\pm$ lie to the left/right of $\Gamma$.

The singular integral operator $S$ with the Cauchy kernel is
defined by
\begin{equation}\label{S} (S\phi)(t)=\frac{1}{\pi i}\oint_\Gamma
\phi(\tau)\frac{d\tau}{\tau -t}. \end{equation}  It acts as an
involution \cite{GKru92} on the linear manifold of all rational
functions with the poles off $\Gamma$, dense in the Hilbert space
${\mathcal H}=L^2(\Gamma)$, with respect to the Lebesque measure
on $\Gamma$. This operator is bounded in $L^2$ norm, and can
therefore be continued to the involution acting on the whole
$L^2(\Gamma)$, if and only if $\Gamma$ is the so called {\em
Carleson curve}. This result, along with the definition of
Carleson curves, as well as detailed proofs and the history of the
subject, can be found in \cite{BK97}. For our purposes it suffices
to know that $S$ is a bounded involution when the curve $\Gamma$
is piecewise smooth, i.e., admits a piecewise continuously
differentiable parametrization.

If $\Gamma$ is a circle or a line, then $S$ is in fact
selfadjoint, and both its norm and essential norm are equal to 1.
This situation is trivial from our point of view, since $W(S)$ and
$W_{\rm ess}(S)$ then coincide with the closed interval $[-1,1]$
and $W_c(S)$ is $[-\norm{c}, \norm{c}]$.

As it happens \cite{Kru91}, circles and lines are the only simple
closed curves in $\dot{\C}$ for which $S$ is selfadjoint. On the
other hand, for all smooth simple closed curves the essential norm
of $S$ is the same, that is, equal to 1 (see \cite[Chapter
7]{GKru92} for Lyapunov curves; the validity of the result for
general smooth curves rests on the compactness result from
\cite{Gru80} and is well known within singular integral
community). Thus, lines and circles are the only smooth closed
curves in $\dot{\C}$ for which the norm and the essential norm of
$S$ coincide. However, such a coincidence is possible for other
{\sl piecewise} smooth (even simple) curves.

One such case occurs when $\Gamma$ is a bundle of $m$ lines
passing through a common point, or of $m$ circles  passing through
two common points. According to \cite{GalKru}, then \[
\norm{S}=\norm{S}_{\rm ess}\geq \cot\frac{\pi}{4m}, \] with the
last inequality turning into equality for at least $m = 1,2,3$.
Respectively, for such curves $\Gamma$ the sets $W(S)$, $W_{\rm
ess}(S)$ are the ellipses with the foci at $\pm 1$, coinciding up
to the boundary, and with the major axes of the length at least
$2\csc\frac{\pi}{2m}$. This length equals $2\csc\frac{\pi}{2m}$
for $m=2,3$. The $c$-numerical range of $S$ is the same ellipse,
only scaled by $\norm{c}$.

The equality $\norm{S}=\norm{S}_{\rm ess}$ also holds for $\Gamma$
consisting of circular arcs (one of which can degenerate into a
line segment) connecting the same two points in $\C$
\cite{Ave,AveKru}; in order for an appropriate orientation on
$\Gamma$ to exist the number of these arcs must be even. If, in
particular, there are two of them (that is, the curve $\Gamma$ is
simple), then \[ \norm{S}=\norm{S}_{\rm ess}= D_\phi
+\sqrt{D^2_\phi+1},\] where \[ D_\phi = \sup \left\{ \frac{\sinh
(\pi\phi\xi)}{\cosh (\pi\xi)}\colon \xi\geq 0\right\} \] and $\pi
(1-\phi)$ is the angle between the arcs forming $\Gamma$
\cite{Ave}. The ellipses $W(S)$, $W_{\rm ess}(S)$ therefore have
the major axes of the length $2\sqrt{D_\phi^2+1}$.

For some particular values of $\phi$ the explicit value of
$D_\phi$ can be easily computed, see \cite{Ave}. If, for instance,
$\Gamma$ consists of a half circle and its diameter, that is $\phi
= 1/2$, then $D_\phi=1/2\sqrt{2}$. Respectively, the major axes of
$W(S)$ and $W_{\rm ess}(S)$ have the length $3/\sqrt{2}$.

It would be interesting to describe all curves $\Gamma$ for which
the norm and the essential norm of the operator (\ref{S}) are the
same.

\subsection{Singular integral operators on weighted spaces on the circle}
Let now $\Gamma$ be the unit circle $\mathbb T$. We
again consider the involution (\ref{S}), this time with ${\mathcal
H}$ being the {\em weighted} Lebesgue space $L^2_{\rho}$. The norm
on this space is defined by \[ \norm{f}_{L^2_{\rho}}=\norm{\rho
f}_{L^2}:=\frac{1}{\sqrt{2\pi}} \left(\int_0^{2\pi}
|f(e^{i\theta})|^2(\rho(e^{i\theta})^2d \theta \right)^{1/2},
\] where the weight $\rho$ is an a.e. positive measurable and square integrable
function on $\mathbb T$. In this setting, the operator $S$ is
closely related with the Toeplitz and Hankel operators on Hardy
spaces, weighted or not. All needed definitions and ``named''
results used below and not supplied with explicit references
conveniently can be found in the exhaustive recent monograph
\cite{Peller}.

\subsubsection{}
Involution $S$ is bounded on $L^2_{\rho}$ if and only if $\rho^2$
satisfies the Helson-Szeg\H{o} condition, that is, can be
represented as
\begin{equation}\label{hs}
\exp(\xi+\overline{\eta}) \text{ with } \xi,\eta\in
L^\infty(\mathbb T) \text{ real valued and }
\norm{\eta}_\infty<\pi/2\end{equation} \cite[p. 419]{Peller}. This
condition is equivalent to
\begin{equation}\label{feb31} \norm{H_\omega}<1,\end{equation} where
\begin{equation}\label{oh} \omega=\overline{\rho_+}/{\rho_+},
\end{equation}$\rho_+$ is the outer
function such that $\abs{\rho_+}=\rho$ a.e. on $\mathbb T$, and
$H_\omega$ denotes the Hankel operator $H_\omega$ with the symbol
$\omega$ acting from the (unweighted) Hardy space $H^2$ to its
orthogonal complement in $L^2$. It is also equivalent to
invertibility of the Toeplitz operator $T_\omega$ on $H^2$.
Moreover \cite{FKS96},
\[
\norm{S}_{L^2_\rho}=\sqrt{\frac{1+\norm{H_\omega}}{1-\norm{H_\omega}}},
\] and a similar relation holds for the essential norms of $S$ and
$H_\omega$. But \[ \norm{H_\omega} = \dist (\omega, H^\infty) \]
(Nehari theorem \cite[p. 3]{Peller}) and  \[ \norm{H_\omega}_{\rm
ess} = \dist (\omega, H^\infty+C) \] (Adamyan-Arov-Krein theorem
\cite[Theorem 1.5.3]{Peller}), where $H^\infty$ is the Hardy class
of bounded analytic in $\mathbb D$ functions, and its sum with the
set $C$ of continuous on $\mathbb T$ functions is the {\em Douglas
algebra} $H^\infty+C$. Thus, the ellipses $W(S)$ and $W_{\rm
ess}(S)$ have the major axes
\[ 2/\sqrt{1-\dist (\omega, H^\infty)} \text{ and }2/\sqrt{1-\dist
(\omega, H^\infty+C)},\] respectively.

The norm of $S$ is attained only simultaneously with the norm of
$H_\omega$. This happens, in particular, if $H_\omega$ is compact,
that is $\omega\in H^\infty+C$. The latter condition can be
restated directly in terms of $\rho$ \cite{FKS96} and means that
$\log\rho\in VMO$, where $VMO$ (the class of functions with
vanishing mean oscillation) is the sum of $C$ with its harmonic
conjugate $\widetilde{C}$.

Thus, for all the weights $\rho$ such that $\log\rho\in VMO$ the
ellipse $W(S)$ is closed, while $W_{\rm ess}(S)$ degenerates into
the line interval $[-1,1]$.

A criterion for the norm of $H_\omega$ to be attained also can be
given, though in less explicit form. Recall that the distance from
$\omega$ to $H^\infty$ is always attained on some $g\in H^\infty$
(this is part of Nehari theorem). This $g$ in general is not
unique, and any $f$ of the form $\omega-g$ is called a {\em
minifunction}. By (another) Adamyan-Arov-Krein's theorem
\cite[Theorem 1.1.4]{Peller}, the norm of $H_\omega$ is attained
if and only if the minifunction is unique and can be represented
in the form
\begin{equation}\label{mini}
f(z)=\norm{H_\omega}\overline{z\theta h}/h,\end{equation} where
$\theta$ and $h\, (\in H^2)$ are some inner and outer functions of
$z$, respectively \footnote{Formally speaking, Theorem 1.1.4 in
\cite{Peller} contains only the ``only if'' part. The ``if''
direction is trivial, since the norm of $H_\omega$ is attained on
$h$ from (\ref{mini}); see Theorem 2.1 of the original paper
\cite{AAK681}.}.

\subsubsection{} We now turn to possible realizations of the
outlined possibilities. If $f$ admits a representation
(\ref{mini}) with $\theta$ of an infinite degree (that is, being
an infinite Blaschke product or containing a non-trivial singular
factor), then $\norm{H_\omega}$ is the $s$-number of $H_\omega$
having infinite multiplicity. In particular,
\begin{equation}\label{eqh} \norm{H_\omega}=\norm{H_\omega}_{\rm
ess}. \end{equation} According to Theorem~\ref{th:wc}, $W(S)$ in
this case coincides with the closed ellipse $W_{\rm ess}(S)$, all
$c$-numerical ranges also are closed and differ from $W(S)$ only
by an appropriate scaling.

Now let $\theta$ in (\ref{mini}) be a finite Blaschke product of
degree $b\,(\geq 0)$ while $h$ is invertible in $H^2$. Suppose
also that $\abs{h}^2$ does not satisfy Helson-Szeg\H{o} condition,
that is, cannot be represented in the form (\ref{hs}) (such outer
functions are easy to construct -- take for example $h$ with
$\abs{h}^{\pm 1}\in L^2$ but $\abs{h}\notin L^{2+\epsilon}$ for
any $\epsilon>0$). Then the Toeplitz operator $T_f$ has
$(b+1)$-dimensional kernel, dense (but not closed) range
\cite[Corollary 3.1 and Theorem 3.16]{LS}, and therefore is not
left Fredholm. By Douglas-Sarason theorem \cite[Theorem 1.1.15]{Peller}
\[ \dist (f, H^\infty+C)=\abs{f}=\norm{H_\omega}=\norm{H_f}.\] We
conclude that (\ref{eqh}) holds again. So, the ellipse $W(S)$ is
closed and coincides with $W_{\rm ess}(S)$. According to
Theorem~\ref{th:wc},  the $c$-numerical range of $S$ is closed if
the number of coefficients $c_j$ of the same sign does not exceed
$b+1$, and open otherwise.

Finally, if a unimodular function $\omega$ is such that the
operator $T_\omega$ is invertible, (\ref{eqh}) holds, but its
minifunction is not constant a.e. in absolute value, then the norm
of $H_\omega$ is not attained. Accordingly, all $c$-numerical
ranges, $W(S)$ in particular, in this case are open.

A concrete realization of the latter possibility is given in the
next subsection. All the other possibilities mentioned earlier
also occur. To construct the respective weights $\rho$, the
following procedure can be applied. Starting with any inner
function $\theta$ and outer function $h\in H^2$, choose $f$ as in
(\ref{mini}) with $\norm{H_\omega}$ changed to an arbitrary
constant in $(0,1)$. Let $\omega$ be an 1-canonical function
\footnote{See \cite[p. 156]{Peller} for the definition.} of the
Nehari problem corresponding to the Hankel operator $H_f$. As
such, $\omega$ is unimodular, and can be represented as
$\omega=g/\overline{g}$, where $g$ is an outer function in $H^2$
\cite[Theorem 5.1.8]{Peller}. Since $\norm{H_\omega}<1$, the
Toeplitz operator $T_{\omega^{-1}}$ is invertible \cite[Theorem
5.1.10]{Peller} (the last two cited theorems from \cite{Peller}
are again by Adamyan-Arov-Krein \cite{AAK681}). The desired weight
is given by $\rho=\abs{g}$.

By Treil's theorem \cite[Theorem 12.8.1]{Peller}, any positive
semi-definite noninvertible operator with zero or infinite
dimensional kernel is unitarily similar to the modulus of a Hankel
operator. Thus, the multiplicity of the norm of $H_\omega$ as its
singular value can indeed assume any prescribed value, whether or
not (\ref{eqh}) holds.

\subsubsection{}
Consider the concrete case of {\em power weights}
\begin{equation}\label{pw} \rho(t)=\prod \abs{t-t_j}^{\beta_j}, \quad t_j\in
{\mathbb T}, \ \beta_j\in\R\setminus \{0\}.\end{equation} It is an
old and well known result that $S$ is bounded on $L^2_\rho$ with
$\rho$ given by (\ref{pw}) if and only if $\abs{\beta_j}<1/2$.
This fact, along with other results about such weights cited and
used below (and established by Krupnik-Verbitskii \cite{VerKru80})
can be found in the monograph \cite[Section 5]{Kru87}.

The essential norm of $S$ does not depend on the distribution of
the nodes $t_j$ along $\mathbb T$, and equals
\begin{equation}\label{sess} \norm{S}_{\rm
ess}=\cot\frac{\pi(1-2\tilde{\beta})}{4}, \text{ where }
\tilde{\beta}=\max\abs{\beta_j}.\end{equation}

In case of only one node (say $t_0$, with the corresponding
exponent $\beta_0$), the norm of $S$ is the same as (\ref{sess}).
The function $\omega$ constructed by this weight $\rho$ in
accordance with (\ref{oh}) is simply $\omega(t)=t^{\beta_0}$,
having a discontinuity at $t_0$. The distance from $\omega$ to
$H^\infty$ is the same as to $H^\infty+C$, it equals $\sin
(\pi\abs{\beta_0})$ and is attained on a constant $\ell=\cos
(\pi\abs{\beta_0})e^{i\pi\beta_0}$. A corresponding minifunction
$f=\omega-\ell$ is not constant a.e. in absolute value; thus, it
cannot admit representation (\ref{mini}). Consequently, the norm
of $H_\omega$ is not attained. Accordingly, $W_c(S)$ is open for
all $c$; the numerical range $W(S)$ has the major axis of the
length $2\sec (\pi\abs{\beta_0})$. Other $c$-numerical ranges are
scaled by $\norm{c}$, as usual.

More generally, the norm of $S$ coincides with (\ref{sess})
independently on the number of nodes, provided that one of the
exponents (say $\beta_0$) differs by its sign from all others and
at the same time exceeds or equals their sum by absolute value.
The size and the shape of all the ellipses $W(S)$, $W_{\rm ess}$,
$W_c(S)$ is then the same as for the weight with only one exponent
$\beta_0$.

In case of two nodes ($t_1$ and $t_2$), the condition above holds
if the respective exponents $\beta_1$, $\beta_2$ are of the
opposite sign. If the signs are the same, the norm of $S$ actually
depends on $\arg t_1/t_2$. It takes its minimal value (for fixed
$\beta_j$) when  $t_1/t_2<0$. This value coincides with
(\ref{sess}), thus making Theorem~\ref{th:wc} applicable again.

\subsection{Composition operators}
For an analytic mapping of the unit disc $\mathbb D$ into itself,
the {\em composition operator} $C_\phi$ is defined as
\[ (C_\phi f)(z)= f(\phi(z)).\]  \subsubsection{} We consider this operator first
on the Hardy space $H^2$. In this setting, the operator $C_\phi$
is bounded and, if $\phi$ is an inner function,
\begin{equation}\label{cn}
\norm{C_\phi}=\sqrt{\frac{1+\abs{\phi(0)}}{1-\abs{\phi(0)}}},
\end{equation} see \cite{Nor68}, also \cite{CoClu95}. It is easily
seen from the proof of (\ref{cn}) given there that the norm of
$C_\phi$ is not attained, unless $\phi(0)=0$. As was shown in
\cite{Sha87, Sha00}, the essential norm of $C_\phi$ coincides with
its norm; moreover, this property is characteristic for inner
functions.

The numerical ranges of composition operators $C_\phi$ with $\phi$
being conformal automorphisms of $\mathbb D$ where treated in
\cite{BoSha00}. It was observed there, in particular, that
$W(C_\phi)$ is an elliptical disc with the foci at $\pm 1$ when
$C_\phi$ is an involution, that is,
\begin{equation}\label{jan311} \phi(z)=\frac{p-z}{1-\overline{p}z}
\end{equation} for some fixed $p\in\mathbb D$. The major axis of this disc $E_p$
was computed in \cite{Abd}, where as a result of rather lengthy
computations it was shown to equal $2/\sqrt{1-\abs{p}^2}$. For
$p=0$, $C_\phi$ is an involution of norm 1. Respectively, $E_0$
degenerates into the closed interval $[-1,1]$. The question of
openness or closedness of $E_p$ for $p\neq 0$ was not discussed.

It follows from Theorem~\ref{th:w} that $E_p$ is open (if $p\neq
0$); moreover, the length of its axes can be immediately seen from
(\ref{invaxes}) and (\ref{cn}): \[
\sqrt{\frac{1+\abs{p}}{1-\abs{p}}}+\sqrt{\frac{1-\abs{p}}{1+\abs{p}}}=
{2}/{\sqrt{1-\abs{p}^2}}.\] Furthermore, Theorem~\ref{th:wess}
implies that $W_{\rm ess}(C_\phi)$ is the closure of $E_p$.
Finally, by Theorem~\ref{th:wc} the $c$-numerical range of
$C_\phi$ is $E_p$ dilated by $\norm{c}$.

\subsubsection{}
These results, with some natural modifications, extend to the case
of {\sl weighted} spaces $H^2_\rho$. Namely, for a non-negative
function $\rho\in L^2({\mathbb T})$ with $\log\rho\in L^1$ we
define the outer function $\rho_+$ as in (\ref{oh}). Then \[
H^2_\rho =\{ f\colon  \rho_+ f\in H^2 \} \text{ and }
\norm{f}_{H^2_\rho}^2=\norm{\rho_+f}_{H^2}.\] A change-of-variable
argument, similar to that used in \cite{Nor68}, shows the
following equality:
\begin{multline} \label{intcomp} \|C_{\phi}f\|_{H^2_\rho}^2=
\frac{1}{2\pi} \int_0^{2\pi}
|f\left(\phi(e^{i\theta})\right)|^2(\rho(e^{i\theta}))^2d \theta \\
=  \frac{1}{2\pi} \int_0^{2\pi}
|f(e^{i\mu})|^2\left(\rho\left(\phi\left(e^{i\mu}\right)\right)\right)^2
\frac{1-|p|^2}{|p-e^{i\mu}|^2}d\mu =\norm{f\chi}^2_{H^2_\rho},
\end{multline}  where \[ \chi(t):=
\frac{\sqrt{1-\abs{p}^2}}
{\abs{p-t}}\frac{\rho(\phi(t))}{\rho(t)}, \quad t\in {\mathbb
T}.\] The norm of a multiplication operator on weighted and
unweighted Hardy spaces is the same. According to (\ref{intcomp})
the operator $C_\phi$ is therefore bounded on $H^2_\rho$ if and
only if
\begin{equation}\label{boundcomp}
\sup_{t\in\mathbb
T}\frac{\rho(\phi(t))}{\rho(t)}<\infty.\end{equation} Observe that
(\ref{boundcomp}) is equivalent to \[ \inf_{t\in\mathbb
T}\frac{\rho(\phi(t))}{\rho(t)}>0 \] because $\phi$ is an
involution. Apparently, (\ref{boundcomp}) holds if $\rho\in
L^\infty$ is bounded below from 0, but there are plenty of
unbounded weights $\rho$ satisfying (\ref{boundcomp}) as well.

Under this condition, $\norm{C_\phi}_{H^2_\rho}=M$, where
\begin{equation}\label{normcphi}
M=\sqrt{1-\abs{p}^2}\sup_{t\in\mathbb
T}\frac{\rho(\phi(t))}{\abs{p-t}\rho(t)}.\end{equation} For any
$\epsilon>0$, consider a function $g\in H^2_\rho$ with the norm
one and such that $\norm{C_\phi g}_{H^2_\rho}>M-\epsilon$. Then
$\norm{C_\phi g_n}_{H^2_\rho}>M-\epsilon$ for $g_n(z)=z^n g(z)$,
$n=1,2,\ldots$. Since the sequence $g_n$ converges weakly to zero
in $H^2_\rho$, from here it follows that the essential norm of
$C_\phi$ also equals $M$.
(We use here the well-known fact that
compact operators on Hilbert spaces
map weakly convergent sequences into strongly convergent
sequences, see \cite[Section 85]{Riesz}, for example.) Moreover, the norm of
$C_\phi$ is attained if and only
if there exist non-zero functions in $H^2_\rho$ with absolute
value equal zero a.e. on the subset of $\mathbb T$ where
$\abs{\chi(t)}\neq M$. Due to uniqueness theorem for analytic
functions, a necessary and sufficient condition for this to happen
is \begin{equation}\label{att}
\abs{\frac{\rho(\phi(t))}{(p-t)\rho(t)}}=\text{const\ a.e.\ on
}\mathbb T.\end{equation} If (\ref{att}) holds, then the norm is
attained in particular on all inner functions, so that the
respective subspace is infinitely dimensional. Consequently,
$W_{\rm ess}(C_\phi)$ is the closed ellipse with the foci at $\pm
1$ and the axes $M\pm M^{-1}$, and $W(C_\phi)$ is the same ellipse
when (\ref{att}) holds or its interior when it does not. The
$c$-numerical range is simply $\norm{c}W(C_\phi)$.

Of course, for $\rho(t)\equiv t$ condition (\ref{boundcomp})
holds, formula (\ref{normcphi}) turns into (\ref{cn}), and
(\ref{att}) is equivalent to $p=0$. Thus, the results obtained
match those already known in the unweighted setting.

\subsubsection{}
One can also consider composition operators $C_\phi$ on weighted
{\sl Lebesgue} spaces $L^2_\rho$. Formula for the norm and the
essential norm of $C_\phi$ remain exactly the same, with no
changes in their derivation\footnote{Moreover, condition
$\log\rho\in L^1$ can be weakened simply to $\rho$ being positive
a.e. on $\mathbb T$, as was the case in Subsection 3.2.}.
Condition for the norm to be attained is different: in place of
(\ref{att}) it is required that the supremum in its left hand side
is attained on a set of positive measure. The respective changes
in the statement about the numerical ranges are evident, and we
skip them. We note only that for $\rho(t)\equiv t$ the supremum in
the right hand side of (\ref{att}) either is attained everywhere
(if $p=0$) or just at one point (if $p\neq 0$). Thus, all the sets
$W(C_\phi)$, $W_{\rm ess}(C_\phi)$ and $W_c(C_\phi)$ are exactly
the same whether the composition operator $C_\phi$ with the symbol
(\ref{jan311}) acts on $H^2$ or $L^2$.

\subsubsection{}
Finally, we consider the operator $C_\phi$ on the Dirichlet space
$\mathcal{D}$. Recall that the latter is defined as the set of all
analytic functions $f$ on $\mathbb D$ such that
$$ \|f\|^2_{\mathcal{D}}:=|f(0)|^2+ \int_{\mathbb D} |f'(z)|^2dA(z) <\infty, $$
where $dA$ is the area measure.

It was shown in \cite[Theorem 2]{MV05} that for any univalent
mapping $\phi$ of $\mathbb D$ onto its subset of full measure,
\[ \|C_{\phi}\|_{\mathcal{D}}=\sqrt{\frac{L+2+\sqrt{L(4+L)}}{2}},
\]
where $L=-\log (1-|\phi(0)|^2)$. This simplifies to \[
\|C_{\phi}\|_{\mathcal{D}}=\frac{\sqrt{L}+\sqrt{4+L}}{2},
\] and is of course applicable when $\phi$ is given by
(\ref{jan311}). Consequently, the elliptical disc  $W(C_{\phi})$
has the major axis \[ \sqrt{4+\log\frac{1}{1-\abs{p}^2}}.\]
Moreover, the operators considered in \cite[Theorem 2]{MV05}
attain their norms, so that $W(C_{\phi})$ is closed.

It was further observed in \cite[Proposition 2.4]{H05} that the
essential norm of $C_\phi$ on ${\mathcal{D}}$ does not exceed $1$,
for any univalent $\phi$. For $\phi$ given by (\ref{jan311}), the
essential norm of $C_\phi$ on ${\mathcal{D}}$ must be equal $1$,
since the essential norm of an involution on an infinite
dimensional space is at least one. Thus, $W_{\rm ess}(C_{\phi})$
in this setting is the closed interval $[-1,1]$.

Analogous remarks can be made in other contexts where the norms
and essential norms of composition operators are known.
\bigskip

{\bf Acknowledgment}. We thank V. Bolotnikov for helpful
discussions concerning composition operators.


\begin{thebibliography}{10}

\bibitem{Abd}
A.~Abdollahi, \emph{The numerical range of a composition operator
with conformal automorphism symbol}, Linear Algebra Appl.
\textbf{408} (2005), 177--188.

\bibitem{AAK681}
V.~M. Adamjan, D.~Z. Arov, and M.~G. Krein, \emph{Infinite
{H}ankel matrices and generalized problems of
{C}arath\'eodory-{F}ej\'er and {F}. {R}iesz}, Funkcional. Anal. i
Prilozhen. \textbf{2} (1968), no.~1, 1--19 (in Russian), English
translation: {\it J. Funct. Anal. and Appl.} {\bf 2} (1968) 1--18.

\bibitem{Ave}
R.~E. Avendan\~o, \emph{Norm and essential norm estimates of
singular integral operators}, Ph.D. thesis, Kishinev State
University, 1988, 109 pp. (in Russian).

\bibitem{AveKru}
R.~E. Avendan\~o and N.~Ya. Krupnik, \emph{A local principle for
calculating quotient norms of singular integral operators},
Funktsional. Anal. i Prilozhen. \textbf{22} (1988), no.~2, 57--58
(in Russian), English translation: {\it Funkt. Anal. Appl.} {\bf
22} (1968), 130--131.

\bibitem{BoDu73}
F.~F. Bonsall and J.~Duncan, \emph{Numerical ranges. {I}{I}},
Cambridge University Press, New York, 1973, London Mathematical
Society Lecture Notes Series, No. 10.

\bibitem{BK97}
A.~B{\"o}ttcher and Yu.~I. Karlovich, \emph{Carleson curves,
{M}uckenhoupt weights, and {T}oeplitz operators}, Birkh\"auser
Verlag, Basel and Boston, 1997.

\bibitem{BoSha00}
P.~S. Bourdon and J.~H. Shapiro, \emph{The numerical ranges of
automorphic composition operators}, J. Math. Anal. Appl.
\textbf{251} (2000), no.~2, 839--854.

\bibitem{ChienTsoWu}
M.-T. Chien, S.-H. Tso, and P.~Y. Wu, \emph{Higher-dimensional
numerical ranges of quadratic operators}, J. Operator Theory
\textbf{49} (2003), no.~1, 153--171.

\bibitem{CoClu95}
C.~C. Cowen and B.~D. MacCluer, \emph{Composition operators on
spaces of analytic functions}, Studies in Advanced Mathematics,
CRC Press, Boca Raton, FL, 1995.

\bibitem{FKS96}
I.~Feldman, N.~Krupnik, and I.~M. Spitkovsky, \emph{Norms of the
singular integral operator with {C}auchy kernel along certain
contours}, Integral Equations and Operator Theory \textbf{24}
(1996), 68--80.

\bibitem{FilStWil72}
P.~A. Fillmore, J.~G. Stampfli, and J.~P. Williams, \emph{On the
essential numerical range, the essential spectrum, and a problem
of {H}almos}, Acta Sci. Math. (Szeged) \textbf{33} (1972),
179--192.

\bibitem{GalKru}
J.~Galperin and N.~Krupnik, \emph{On the norms of singular
integral operators along certain curves with intersections},
Integral Equations and Operator Theory \textbf{29} (1997), no.~1,
10--16.

\bibitem{GKru92}
I.~Gohberg and N.~Krupnik, \emph{One-dimensional linear singular
integral equations. {I}ntroduction}, OT 53, 54, vol. 1 and 2,
Birkh{\"a}user Verlag, Basel and Boston, 1992.

\bibitem{Gru80}
S.~M. Grudsky, \emph{On the compactness of a certain integral
operator}, No. 4856-80 dep., VINITI, Moscow, 1980 (in Russian).

\bibitem{GusRa}
K.~E. Gustafson and D.~K.~M. Rao, \emph{Numerical range. {T}he
field of values of linear operators and matrices}, Springer, New
York, 1997.

\bibitem{Ha}
P.~R. Halmos, \emph{A {H}ilbert space problem book}, Van Nostrand,
Princeton, NJ, 1967.

\bibitem{H05}
C.~Hammond, \emph{The norm of a composition operator with linear
symbol acting on the dirichlet space}, J. Math. Anal. Appl.
\textbf{303} (2005), 499--508.

\bibitem{Kru91}
N.~Krupnik, \emph{The conditions of selfadjointness of the
operator of singular integration}, Integral Equations and Operator
Theory \textbf{14} (1991), 760--763.

\bibitem{Kru87}
N.~Ya. Krupnik, \emph{Banach algebras with symbol and singular
integral operators}, Birkh{\"{a}}user, Basel and Boston, 1987.

\bibitem{LS}
G.~S. Litvinchuk and I.~M. Spitkovsky, \emph{Factorization of
measurable matrix functions}, OT25, Birkh{\"{a}}user Verlag, Basel
and Boston, 1987.

\bibitem{MV05}
M.~J. Mart{\i}n and D.~Vukoti{\'{c}}, \emph{Norms and spectral
radii of composition operators acting on the dirichlet space}, J.
Math. Anal. Appl. \textbf{304} (2005), 22--32.

\bibitem{Nor68}
E.~A. Nordgren, \emph{Composition operators}, Canad. J. Math.
\textbf{20} (1968), 442--449.

\bibitem{Peller}
V.~V. Peller, \emph{Hankel operators and their applications},
Springer, New York-Berlin-Heidelberg, 2003.

\bibitem{Riesz}
F.~Riesz and B.~Sz.-Nagy, \emph{Functional analysis}, Frederick
Ungar Publishing Co., New York, 1955.

\bibitem{Sha87}
J.~H. Shapiro, \emph{The essential norm of a composition
operator}, Annals of Math. \textbf{125} (1987), 375--404.

\bibitem{Sha00}
\bysame, \emph{What do composition operators know about inner
functions?}, Monatsh. Math. \textbf{130} (2000), no.~1, 57--70.

\bibitem{Spit80}
I.~M. Spitkovsky, \emph{Some estimates for partial indices of
measurable matrix valued functions}, Mat. Sb. (N.S.)
\textbf{111(153)} (1980), no.~2, 227--248, 319 (in Russian),
English translation: {\it Math. USSR Sbornik} {\bf 39}
  (1981), 207--226.

\bibitem{Spitkov}
\bysame, \emph{Once more on algebras generated by two
projections}, Linear Algebra Appl. \textbf{208/209} (1994),
377--395.

\bibitem{StWil68}
J.~G. Stampfli and J.~P. Williams, \emph{Growth conditions and the
numerical range in a {B}anach algebra}, T\^ohoku Math. J. (2)
\textbf{20} (1968), 417--424.

\bibitem{TsoWu}
S.-H. Tso and P.~Y. Wu, \emph{Matricial ranges of quadratic
operators}, Rocky Mountain J. Math. \textbf{29} (1999), no.~3,
1139--1152.

\bibitem{VerKru80}
I.~E. Verbickii and N.~Ya. Krupnik, \emph{Exact constants in
theorems on the boundedness of singular operators in {$L\sb{p}$}
spaces with a weight and their application}, Mat. Issled. (1980),
no.~54, 21--35, 165 (in Russian).

\bibitem{West75}
R.~Westwick, \emph{A theorem on numerical range}, Linear and
Multilinear Algebra \textbf{2} (1975), 311--315.

\end{thebibliography}
\providecommand{\bysame}{\leavevmode\hbox
to3em{\hrulefill}\thinspace}
\providecommand{\MR}{\relax\ifhmode\unskip\space\fi MR }
\providecommand{\MRhref}[2]{%
  \href{http://www.ams.org/mathscinet-getitem?mr=#1}{#2}
} \providecommand{\href}[2]{#2}

\end{document}